\documentclass[10pt]{amsart}
\usepackage{amsmath,amsthm,amssymb, amscd, amsfonts}

\usepackage[dvipsnames]{xcolor}

\theoremstyle{plain}
\newtheorem{theorem}{Theorem}[section]
\newtheorem{corollary}[theorem]{Corollary}

\newtheorem{proposition}[theorem]{Proposition}
\newtheorem{lemma}[theorem]{Lemma}
\newtheorem{question}[theorem]{Question}
\theoremstyle{definition}
\newtheorem{definition}[theorem]{Definition}

\newtheorem{example}[theorem]{Example}

\theoremstyle{remark}
\newtheorem{remark}[theorem]{Remark}

\numberwithin{equation}{section}\theoremstyle{plain}

\newcommand{\C}{{\mathbb C}}
\newcommand{\A}{{\mathcal A}}

\newcommand\id{\operatorname{id}}
\newcommand\dd{\operatorname{d}}
\newcommand\Supp{\operatorname{Supp}}

\newcommand\co{\operatorname{co}}

\newcommand\cop{\operatorname{cop}}

\newcommand\Opext{\operatorname{Opext}}

\newcommand\op{\operatorname{op}}

\begin{document}
\title[Finite quantum groups and quantum permutation groups]{Finite quantum groups and quantum permutation groups}
\author{Teodor Banica}
\address{T. Banica: D\' epartement de  Math\'ematiques, Universit\'e Cergy-Pontoise, 95000 Cergy-Pontoise,
France.} \email{teodor.banica@u-cergy.fr}
\author{Julien Bichon}
\address{J. Bichon:
Laboratoire de Math\'ematiques, Universit\'e Blaise Pascal, 63177
Aubi\`ere Cedex, France.} \email{
Julien.Bichon@math.univ-bpclermont.fr }
\author{Sonia Natale}
\address{S. Natale: Fa.M.A.F.
Universidad Nacional de C\'ordoba. CIEM -- CONICET. (5000)
C\'ordoba, Argentina.} \email{natale@famaf.unc.edu.ar}

\thanks{Partially supported by CONICET--CNRS,  SeCYT (UNC), and the ANR project Galoisint.}


\date{April 7, 2011}

\begin{abstract} We give examples of finite quantum permutation groups which arise from the twisting construction or as bicrossed products  associated
to exact factorizations in finite groups. We also give examples of
finite quantum groups which are not  quantum permutation groups:
one such example occurs as a split abelian extension associated to
the exact factorization $\mathbb S_4 = \mathbb Z_4 \mathbb S_3$ and has
dimension $24$. We show that, in fact, this is the smallest possible dimension that a non quantum permutation group can have.
\end{abstract}

\maketitle

\section{Introduction}




Let $n \geq 1$ be an integer. Recall from Wang's paper \cite{wang} that the usual
symmetric group $\mathbb S_n$ has a free analogue, denoted
$\mathbb S_n^+$; this is a compact quantum group acting universally on the set $\{1, \dots, n\}$.

We shall work over an algebraically closed base field $k$ of
characteristic zero.
Let  $\A_s(n) = \A_s(n, k)$ be the Hopf algebra corresponding to Wang's quantum permutation group
\cite{bichon}. This is the algebra given by generators $u_{ij}$, $1\leq i, j \leq
n$, with relations making $(u_{ij})_{i, j}$ a \emph{magic} matrix, that is,
\begin{equation}\label{qpa1}u_{ij}u_{ik} = \delta_{jk}u_{ij}, \quad u_{ij}u_{kj} = \delta_{ik}u_{ij}, \quad \sum_{l = 1}^nu_{il} = 1
= \sum_{l = 1}^nu_{li},\end{equation} for all $1\leq i, j, k \leq
n$.  The algebra $\A_s(n)$ is a cosemisimple Hopf algebra with comultiplication, counit and antipode
determined by
\begin{equation}\label{qpa2}\Delta (u_{ij})
= \sum_{k = 1}^nu_{ik} \otimes u_{kj}, \quad  \epsilon(u_{ij}) =
\delta_{ij}, \quad \mathcal S(u_{ij}) = u_{ji}, \quad 1 \leq i, j
\leq n. \end{equation}

By a (finite) \emph{quantum permutation algebra} we
shall understand a (finite dimensional) quotient Hopf algebra $H$
of $\A_s(n)$.  In particular, a quantum permutation algebra satisfies $\mathcal S^2 = \id$.
Hence, a finite quantum permutation algebra is cosemisimple and semisimple.
Formally, a quantum permutation algebra corresponds to a ``quantum permutation group'', that is, a quantum subgroup of $\mathbb S_n^+$.

The Hopf algebra $\A_s(n)$ is the
universal cosemisimple Hopf algebra coacting on the commutative
algebra $k^n$ \cite{bichon}. Thus, a cosemisimple Hopf algebra $H$ is a
quantum permutation algebra if and only if there exists a separable commutative  faithful (left or right) $H$-comodule algebra
$L$.

The Hopf algebra $k^{\mathbb S_n}$ of functions on the classical symmetric group is a quantum permutation algebra (the maximal commutative quotient of $\A_s(n)$). We have $k^{\mathbb S_n} = \A_s(n)$, for $n = 1, 2, 3$, but not for $n \geq 4$, where the algebra $\A_s(n)$ is not commutative and infinite dimensional.
See \cite{wang, bichon}.




Several interesting examples of quantum permutation
groups, finite or not, were obtained as twisting deformations of classical groups. Here is a list
of such quantum groups:

(1) The twists of $\mathbb S_n$ from \cite{bi1}.

(2) The quantum group $O_n^{-1}$, see \cite{ahn}.

(3) The twists of several subgroups of $SO_3$, see \cite{qfp}.

Also the nontrivial Hopf algebras studied by Masuoka \cite{ma-contemp} appear as quantum permutations
algebras in \cite{qfp}, including the historical 8-dimensional Kac-Paljutkin example (which is not a twist
of a function algebra).

\medbreak As is well-known, the Cayley representation makes every finite group into a permutation group.  This leads naturally to consider the question whether
any ``finite quantum group'' is a ``quantum permutation group''. In other words:
\begin{question}\label{q1}
Let $H$ be a finite dimensional cosemisimple Hopf algebra. Is it true that $H$ is a quotient of $\A_s(N)$, for some $N \in \mathbb N$?
\end{question}

Of course, the answer is 'yes' if $H$ is commutative, taking $N = \dim H$.
Although less trivial to see, the answer is also 'yes' in the cocommutative case.
Moreover, as shown in
\cite{bichon}, the cocommutative cosemisimple (finite or not) Hopf algebra quotients of $\A_s(N)$ are exactly the group algebras $kG$,
where $G$ is a quotient of a free product $G_1 * \dots * G_m$ of transitive abelian groups $G_i \subseteq \mathbb S_{n_i}$, with $N=\sum_i n_i$.

Now back to our question, let $G$ be a finite group of order
$n$. Then there is a canonical surjective homomorphism $\mathbb Z_n^{*n}\to G$,
given by $k^{(g)}\mapsto g^k$. Thus $kG$ is a quotient of $\A_s(N)$, with $N = n^2$.

This also suggests that one could include the condition $N = (\dim H)^2$ in Question \ref{q1}.

\medbreak In this paper we show that the answer to Question \ref{q1} is negative in general. More precisely, we give examples of finite dimensional cosemisimple Hopf algebras which are not quantum permutation algebras. Such examples  arise as split abelian extensions  from  exact factorizations of the
symmetric groups $\mathbb S_4$ and $\mathbb S_5$. See Theorem \ref{s5}.

These examples show that the class of finite quantum permutation algebras is not stable under extensions.
It turns out that their duals are quantum permutation algebras, so we  get that the class of finite quantum permutation algebras is also not stable under
duality. An argument involving Drinfeld doubles implies, in addition, that this class is not stable under twisting deformations neither.

\medbreak We also discuss sufficient conditions on  abelian extensions or a twisting deformation of a linear algebraic group in
order that they be quantum permutation algebras. This is done in Sections
\ref{examples-mp} and \ref{twist}, respectively. Some known examples of quantum permutation algebras turn out to fit into these pictures.

We show that central abelian extensions and certain classes of split extensions, that include cocentral split extensions and split extensions by an abelian group, are quantum permutation algebras.  See Theorems \ref{central} and \ref{f-abelian}.

As a consequence, we get that if $G$ is a finite group, then the Drinfeld double $D(G)$ and its dual $D(G)^*$ are quantum permutation algebras. We also obtain that a cosemisimple Hopf algebra whose dimension divides $p^3$ or $pqr$, where $p$, $q$ and $r$ are pairwise distinct prime numbers, is a quantum permutation algebra (Proposition \ref{pqr}). Other known examples also fit into this picture, like, for instance, some nontrivial Hopf algebras studied by Masuoka \cite{ma-contemp}.

\medbreak We then look at twisting deformations of  a quantum permutation algebra  $H$. We give in Proposition \ref{suff-twist} a general sufficient condition on a cocycle $\sigma: H \otimes H \to k$ such that the twisted Hopf algebra $H^\sigma$ is a quantum permutation algebra. The deformations of the symmetric groups  in \cite{bi1} fall into this class.

Let $\Gamma$ be a finite abelian group and $\sigma$ a 2-cocycle on $\Gamma$.  We give further examples of quantum permutation algebras as twisting deformations of linear algebraic groups $G$, with $\widehat{\Gamma} \subset G \subset {\rm Aut}(k_\sigma\Gamma)$. This construction relies on the results of \cite{tga}. See Theorem \ref{twist-abelian}. The twisted examples from \cite{ahn} fit into this framework.

\medbreak In Section \ref{QPE} we introduce the \emph{quantum permutation envelope}, denoted $H_{qp}$, of a cosemisimple Hopf algebra as the  subalgebra generated by the matrix coefficients of all separable commutative (right and left) coideal subalgebras of $H$. This is a Hopf subalgebra containing all quantum permutation algebras $A \subseteq H$. When $H$ is finite dimensional, $H_{qp}$ is the maximal quantum permutation algebra contained in $H$. This provides a method to construct quantum permutation algebras from any finite dimensional cosemisimple Hopf algebra. We determine the quantum permutation envelope for some families of examples.

\subsection*{Conventions} We refer the reader to \cite{Mo} for the notation and terminology on Hopf algebras used throughout. By a \emph{twisting deformation} of a Hopf algebra $H$ we understand a twist $H^\sigma$ in the sense of Doi \cite{doi}, where $\sigma: H \otimes H \to k$ is a convolution invertible normalized $2$-cocycle. That is, $H^\sigma = H$ as a coalgebra, with multiplication
$$[x] [y]= \sigma(x_{1}, y_{1})
\sigma^{-1}(x_{3}, y_{3}) [x_{2} y_{2}], \quad x,y \in H,$$ where $[x]$ denotes the element $x \in H$, viewed as an element of $H^{\sigma}$.

\section{Quantum permutation algebras}\label{QPA}

Let $H$ be a quantum permutation algebra. As noticed before, we have $\mathcal S^2 = \id$ in $H$.
For $H$ finite-dimensional this condition is equivalent to $H$ being separable and/or cosemisimple \cite{LR}.

\begin{definition} The \emph{degree} of $H$, denoted
$\dd (H)$,  is the smallest $n \geq 1$, such that $H$ is a quotient Hopf algebra of $\A_s(n)$. \end{definition}

As pointed out in the Introduction,  if $H = k^G$, where $G$ is a finite group, then
$\dd (H) \leq |G|$, while if $H = kG$, then $\dd (H) \leq |G|^2$ (see
\cite[Proposition 5.2]{bichon}).

\medbreak Let  $f(u_{ij}) = : x_{ij} \in H$, where $f: \A_s(n)
\to H$ is a surjective Hopf algebra map.

For each $j = 1, \dots, n$, consider the subspace $L^j$ of $H$
spanned by $x_{ij}$, $1 \leq i \leq n$. By \eqref{qpa1} and
\eqref{qpa2}, each $L^j$ is a commutative separable left coideal subalgebra
of $H$.  Similarly, the subspace $R^j = \mathcal S(L^j)$ of $H$
spanned by $x_{ji}$, $1 \leq i \leq n$, is a commutative separable right
coideal subalgebra.

In particular, $H$ is generated as an algebra by its commutative separable
left coideal subalgebras $L^1, \dots, L^n$.

\begin{remark}\label{gen-qpg} If $H_1, H_2$ are quantum permutation algebras, then so is their free product
$H_1 * H_2$ (with block-diagonal magic matrix) \cite{wang, repns}.  It follows that any
cosemisimple Hopf algebra $H$ such that $H$ is generated as an
algebra by a finite number of quantum permutation algebras is
itself a quantum permutation algebra.

Further, if $H = k[H_1, \dots, H_r]$ is generated by the quantum
permutation algebras $H_1, \dots, H_r$, then $\dd (H) \leq \dd
(H_1) + \dots + \dd (H_r)$. \end{remark}


\begin{theorem}\label{gen-cmc} Let $H$ be a finite dimensional cosemisimple Hopf algebra.
Then $H$ is a quantum permutation algebra if and only if $H$ is
generated, as an algebra, by the matrix coefficients of its
commutative left (or right) coideal subalgebras.
\end{theorem}

\begin{proof} We have already proved the 'only if' implication.
To prove the converse, suppose $L\subseteq H$ is a left (or right) coideal subalgebra. It is known that $H$ is free as a (left or right) $L$-module under  multiplication (see \cite{ma-freeness} or, more generally, \cite{skryabin}). This implies that $L$ is a separable algebra, by \cite{ma-csub}.

Let $H[L]$ be the subalgebra generated by the matrix coefficients of $L$. Then $H[L]$ is a subbialgebra of $H$, and therefore a Hopf subalgebra, because it is finite dimensional. Suppose $L$ is commutative. Since, by definition,  $H[L]$ coacts faithfully on $L$, then it is a quantum permutation algebra \cite{bichon}.

Since $H$ is finite dimensional, the assumption implies that $H$ is generated by a finite number of the Hopf subalgebras $H[L]$. Hence $H$ is a quantum permutation algebra, by Remark \ref{gen-qpg}. This finishes the proof of the theorem. \end{proof}

If $H$ is a quantum permutation algebra, then so are $H^{\op}$ and $H^{\cop}$.
For a finite dimensional Hopf algebra $H$, the Drinfeld
double $D(H)$ is generated as an algebra by $H^{*\cop}$  and $H$. Moreover, $D(H)$ is cosemisimple if $H$ (and therefore also $H^*$) is cosemisimple.
Then we get:

\begin{corollary}\label{dh} Let $H$ be a finite dimensional
cosemisimple Hopf algebra.  If $H$ and $H^*$ are quantum
permutation algebras, then the Drinfeld double $D(H)$ is a quantum
permutation algebra and we have $\dd (D(H)) \leq \dd (H) + \dd
(H^*)$.

In particular, every Drinfeld double $D(G)$ of a finite
group algebra $kG$ is a quantum permutation algebra  of degree at
most $|G|(1 + |G|)$.  \end{corollary}

We shall see later (Corollary \ref{dual-double}) that $D(G)^*$ is also a quantum permutation algebra.

\section{Hopf algebra extensions}

Recall that an \emph{exact sequence} of finite dimensional Hopf algebras is a sequence of Hopf algebra maps \begin{equation}\label{seq}k \to A \overset{\iota}\to H \overset{\pi}\to \overline H \to k,\end{equation} where $H$ is finite dimensional, such that $\iota$ is injective, $\pi$ is surjective and, identifying $A$ with a Hopf subalgebra of $H$, we have \begin{equation}\label{kernel}A = H^{\co \pi} = \{ h \in H|\, (\id \otimes \pi) \Delta(h) = h \otimes 1 \}.\end{equation}
See \cite[Definition 1.4]{ma-newdir} for details. (In particular,  condition \eqref{kernel} is equivalent to $H/HA^+ = \overline H$, where $A^+ = \ker \epsilon_A$, and we identify $\overline H$ as a quotient Hopf algebra of $H$.) We shall say in this case that $H$ is an \emph{extension} of $A$ by $\overline H$.

\begin{remark}\label{dim-sec} Suppose $H$ is finite dimensional and $A \subseteq H^{\co \pi}$ is a Hopf subalgebra. Let us point out that in this case exactness of the sequence \eqref{seq} is equivalent to the condition $\dim A \dim \overline H = \dim H$. \end{remark}

The following result on Hopf algebra extensions will be used  repeatedly.

\begin{theorem}\label{x-lx} Suppose $H$ is finite dimensional cosemisimple. Assume in addition that:

\emph{(i)} $A$ is a quantum permutation algebra, and

\emph{(ii)} $\overline H$ is generated as an algebra by a finite subset $X$
such that, for all $x \in X$, there exists a commutative left
coideal subalgebra $L^x$ of $H$ with $x \in \pi(L^x)$.

Then $H$ is a quantum permutation algebra. Moreover, we have $\dd (H) \leq \dd(A) + \sum_{1\neq x \in X} \dim
L^x$. \end{theorem}

Note that the same statement holds true replacing \emph{left} by
\emph{right}.

\begin{proof} For each $x \in X$, let $H^x \subseteq H$ be the subalgebra generated by the
subcoalgebra $C^x$ of matrix coefficients of $L^x$. Thus $H^x$ is
a Hopf subalgebra, and we have $L^x \subseteq H^x$.  Since $L^x$
is commutative, then $H^x$ is a quantum permutation algebra, and
$\dd (H^x) \leq \dim L^x$, by construction.

Let $\tilde H$ be the subalgebra generated by  $A$ and $H^x$,
$1\neq x \in X$.  This is a Hopf subalgebra of $H$ that contains
$A$. Moreover, since every element $1\neq x \in X$ belongs to
$\pi(L^x) \subseteq \pi(\tilde H)$, and $X$ generates $\overline
H$, then $\pi\vert_{\tilde H}:\tilde H \to \overline H$ is
surjective.

On the other hand, by exactness of the sequence $k \to A \to H \overset{\pi}\to \overline H \to k$, we have $H^{\co \pi} = A$. Thus $\tilde{H}^{\co \pi|_{\tilde{H}}} = H^{\co \pi} \cap \tilde{H} = A$, since $A \subseteq \tilde{H}$. Therefore the sequence  $k \to A \to \tilde{H} \to \overline H \to k$ is also exact. Then $\tilde H = H$, since they have the same finite  dimension (see Remark \ref{dim-sec}). Thus $H$ is generated as an
algebra by the quantum permutation algebras $A$, $H^x$, $1 \neq x
\in X$. Hence $H$ is a quantum permutation algebra, with the
claimed bound for $\dd (H)$, by Remark \ref{gen-qpg}.
\end{proof}

\section{Matched pairs of groups}

Let $(F, \Gamma)$ be a \emph{matched pair} of finite groups. That
is, $F$ and $\Gamma$ are endowed  with actions by permutations
$\Gamma \overset{\vartriangleleft}\leftarrow \Gamma \times F
\overset{\vartriangleright}\to F$ such that
\begin{equation}\label{matched}
s \vartriangleright xy  = (s \vartriangleright x) ((s
\vartriangleleft x) \vartriangleright y), \quad st
\vartriangleleft x  = (s \vartriangleleft (t \vartriangleright x))
(t \vartriangleleft x), \end{equation} for all $s, t \in \Gamma$,
$x, y \in F$.

Given finite groups $F$ and $\Gamma$, providing them with a pair
of compatible actions is equivalent to giving a group $G$ together
with an \emph{exact factorization} $G = F \Gamma$: the relevant
actions are determined by the relations $gx = (g \vartriangleright
x)(g \vartriangleleft x)$, $x \in F$, $g \in \Gamma$.

\medbreak Consider the left action of $F$ on $k^\Gamma$, $(x.
f)(g) = f(g \vartriangleleft x)$, $f \in k^\Gamma$, and let
$\sigma: F \times F \to (k^*)^\Gamma$ be a normalized
2-cocycle. Dually, consider the right action of $\Gamma$ on
$k^F$, $(w.g)(x) = w(x \vartriangleright g)$, $w \in k^F$, and let
$\tau: \Gamma \times \Gamma \to (k^*)^F$ be a normalized
2-cocycle.

Under appropriate compatibility conditions between $\sigma$ and
$\tau$, the vector space $H = k^\Gamma \otimes k F$  becomes a
(semisimple) Hopf algebra, denoted $H = k^\Gamma \,
{}^{\tau}\#_{\sigma}kF$, with the crossed product algebra
structure and the crossed coproduct coalgebra structure (see \cite[Section 1]{ma-newdir}). For all
$g,h\in \Gamma$, $x, y\in F$, we have
\begin{align}\label{mult} (e_g \# x)(e_h \# y) & = \delta_{g \vartriangleleft x, h}\, \sigma_g(x,
y) e_g \# xy, \\
\label{delta} \Delta(e_g \# x) & = \sum_{st=g} \tau_x(s, t)\, e_s
\# (t \vartriangleright x) \otimes e_{t}\# x,
\end{align}where $\sigma_s(x, y) = \sigma(x, y)(s)$ and $\tau_x(s, t) = \tau(s, t)(x)$, $s, t\in \Gamma$, $x, y\in F$.

Let $\pi = \epsilon \otimes \id: H = k^\Gamma \,
{}^{\tau}\#_{\sigma}kF \to kF$ denote the canonical projection. We
have an exact sequence of Hopf algebras $k \to k^\Gamma \to H
\overset{\pi}\to kF \to k$. Moreover, every Hopf algebra $H$
fitting into an exact sequence of this form is isomorphic to
$k^\Gamma \, {}^{\tau}\#_{\sigma}kF$ for appropriate compatible
actions and cocycles $\sigma$ and $\tau$. Equivalence classes of
such extensions associated to a fixed matched pair  $(F, \Gamma)$
form an abelian group $\Opext(k^\Gamma, kF)$, whose unit element
is the class of the \emph{split} extension $k^\Gamma \# kF$.

\begin{remark} Suppose $k = \C$ is the field of
complex numbers. Then $H$ is a Hopf $C^*$-algebra (often called \emph{Kac algebra}), that is,  it is
a $C^*$-algebra such that all structure maps are $C^*$-algebra
maps.   \cite{kac, ma-ff}. \end{remark}

\begin{remark}\label{h-sub} Let $H = k^\Gamma \, {}^{\tau}\#_{\sigma}kF$ be a bicrossed product.
Let $F' \subseteq F$ be a subgroup, and consider the subgroup $\Gamma' \subseteq \Gamma$ consisting of all elements $g \in \Gamma$ such that $g \vartriangleright F' = F'$. Then $(\Gamma', F')$ is a matched pair by restriction. Indeed, if $s\in \Gamma'$, $x, y \in F'$, then it follows from the compatibility between $\vartriangleright$ and $\vartriangleleft$ that $$(s \vartriangleleft x) \vartriangleright y = (s \vartriangleright x)^{-1} (s \vartriangleright xy) \in F',$$ whence $s \vartriangleleft x \in \Gamma'$. Hence $\Gamma'$ is $F'$-stable, which implies the claim.

In particular, if $F' \subseteq F$ is a subgroup stable under the
action $\vartriangleright$, then $(\Gamma, F')$ is a  matched pair
by restriction, and it follows from formulas \eqref{mult} and
\eqref{delta} that the bicrossed product $k^{\Gamma}\,
{}^{\tau'}\#_{\sigma'}kF'$ is naturally a Hopf subalgebra of $H$,
where $\sigma'$ is the restriction of $\sigma$ to $F'$, and $\tau'_x(g, h) = \tau_x(g, h)$, for all $x\in F'$, $g, h \in
\Gamma$.

Observe that if $F' \subseteq F$ is the largest subgroup acting trivially on
$\Gamma$, then $F'$ is $\Gamma$-stable, by \eqref{matched}. The
Hopf subalgebra $k^{\Gamma}\, {}^{\tau'}\#_{\sigma'}kF'$ is in this
case a central extension. \end{remark}

\section{Quantum permutation algebras obtained from matched pairs of
groups}\label{examples-mp}
Consider a bicrossed product $H = k^\Gamma \,
{}^{\tau}\#_{\sigma}kF$. We shall give in this section sufficient conditions in order for $H$ to be a quantum permutation
algebra. These include the following cases:

\begin{itemize}\item[(1)] $H$ is a central abelian extension (that is, $k^\Gamma$ central in $H$).
\item[(2)] $H$ is a split abelian extension (that is, $\sigma = 1$, $\tau = 1$)
and $F$ is generated by its abelian $\Gamma$-stable subgroups.
(In particular, this is true when $F$ is abelian or the action $\vartriangleright: \Gamma \times F \to F$ is trivial.)
\end{itemize}

The result for Case (1) is a consequence of Theorem \ref{x-lx}:

\begin{theorem}\label{central} Let $H = k^\Gamma \, {}^{\tau}\#_{\sigma}kF$ and suppose that $k^\Gamma$ is central in $H$.
Then $H$ is a quantum permutation algebra and we have $\dd (H)
\leq |\Gamma|\, |F|^2$.  \end{theorem}

\begin{proof} The assumption that $k^\Gamma$ is central implies that the action $\vartriangleleft: \Gamma \times F \to \Gamma$ is trivial.

Let $x \in F$ and let $\langle x\rangle \subseteq F$  denote the
cyclic subgroup generated by $x$. Consider the subspace $L^x =
k^\Gamma \# \langle x\rangle \subseteq H$. It follows from
\eqref{mult} and \eqref{delta} that $L^x$ is a left coideal
subalgebra of $H$. As an algebra, $L^x = k^\Gamma
\#_{\overline\sigma}k\langle x\rangle$ is a crossed product with
respect to the trivial action $k\langle x\rangle \otimes k^\Gamma
\to k^\Gamma$ and the $2$-cocycle $\overline \sigma =
\sigma\vert_{\langle x\rangle \times \langle x\rangle}$. Therefore
$L^x$ is a commutative left coideal subalgebra of dimension
$|x||\Gamma|$.

It is clear that $x \in \pi(L^x)$, for all $x\in F$. Hence $H$ is
a quantum permutation algebra, by Theorem \ref{x-lx}. Moreover,  we
have that $\dd (H) \leq \sum_{x \in F} \dim L_x = |\Gamma| \sum_{x
\in F} |x| \leq  |\Gamma|\, |F|^2$, as claimed.
\end{proof}


Consider next a \emph{split} abelian extension $H = k^\Gamma \# kF$.
It follows from \eqref{mult} and \eqref{delta} that for
any subgroup $T \subseteq F$ such that $T$ is stable under the
action $\vartriangleright$ of $\Gamma$, the group algebra $kT
\simeq 1 \# kT$ is a right coideal  subalgebra of $H$.


\begin{theorem}\label{f-abelian} Let $H = k^\Gamma \# kF$ be a split abelian extension.
Suppose  $F$ is generated by its abelian $\Gamma$-stable subgroups. Then $H$ is a quantum permutation algebra.
Furthermore, we have $\dd (H) \leq |\Gamma|   |F|^2$.  \end{theorem}

\begin{proof} We may assume that $|F|, |\Gamma| > 1$.
Let $T_1, \dots, T_s$, be abelian $\Gamma$-stable subgroups  of
$F$,  with $F = \langle T_1, \dots, T_s\rangle$. Then $kT_i \simeq
1 \# kT_i$ are commutative right coideal subalgebras of $H$, and
$x \in \pi(kT_i)$, for all $x \in T_i$. It follows from  Theorem
\ref{x-lx} that $H$ is a quantum permutation algebra. Moreover, we
have \begin{equation*}\dd (H) \leq |\Gamma| + \sum_{x \in
\cup_iT_i}|T_i^x| \leq |\Gamma| + |F| \, |\cup_i T_i| \leq
|\Gamma| + |F|^2 \leq  |\Gamma| |F|^2, \end{equation*} where in
the first inequality, $T_i^x$ denotes a choice of one the
subgroups $T_i$ such that $x \in T_i$. \end{proof}

\begin{remark} Let $T \subseteq F$ be any abelian subgroup.
Consider the subgroup $\Gamma_T \subseteq \Gamma$ consisting of
all  elements $g \in \Gamma$ such that $g \vartriangleright T =
T$. Then $(\Gamma_T, T)$ is a matched pair by restriction (see
Remark \ref{h-sub}).
Theorem  \ref{f-abelian} implies that the associated split extension $k^{\Gamma_T} \# kT$ is a quantum permutation algebra.
\end{remark}

\begin{remark} Observe that the conclusion in Theorem \ref{f-abelian} holds in either of the following cases:

(i) $F$ is abelian, or

(ii) the action $\vartriangleright: \Gamma \times F \to F$ is trivial, that is, $H$ is a split cocentral extension. \end{remark}


We next discuss some families of examples of finite quantum permutation algebras.

\begin{example} The dual of the Drinfeld double $D(G)$ of a finite group $G$ fits into a central abelian exact sequence $k \to k^G \to D(G)^* \to kG \to k$. By Theorem \ref{central}, we get:

\begin{corollary}\label{dual-double} Let $G$ be a finite group. Then $D(G)^*$ is a quantum permutation algebra and $\dd(D(G)^*) \leq |G|^3$. \end{corollary}
\end{example}

\begin{example} \emph{Dimension $pqr$.} Let $p$, $q$ and $r$ be pairwise distinct prime
numbers. A semisimple Hopf algebra $H$ of
dimension $p$, $p^2$ or $pq$ is necessarily commutative or cocommutative,
so $H$ is a quantum permutation algebra.

\medbreak It is known that every semisimple Hopf algebra $H$ of
dimension $p^3$ fits into a central abelian extension $k \to
k^{\mathbb Z_p} \to H \to k(\mathbb Z_p \times \mathbb Z_p) \to k$
\cite{masuoka-pp}. Hence $H$ is a quantum permutation algebra.

\medbreak  Assume next that $\dim H = pqr$. By \cite[Corollary 9.4]{ENO2}, $H$ is a split abelian extension. Such extensions are classified in \cite[Section 4]{pqq}; in particular, they must be either central or cocentral. It follows from Theorems \ref{central} and \ref{f-abelian} that $H$ is a quantum
permutation algebra.

In conclusion, we can state the following:

\begin{proposition}\label{pqr} Suppose that the dimension of $H$ divides $p^3$ or $pqr$. Then $H$ is a quantum permutation algebra.  \end{proposition}

Consider the case where $\dim H = pq^2$, $p \neq q$. By
the results in \cite[Subsection 9.2]{ENO2} and the classification
results of abelian extensions in \cite{pqq}, either $H$ or $H^*$ fits into a central abelian exact sequence.

\begin{proposition} Suppose that $H$ is nontrivial and $\dim H = pq^2$, $p \neq q$. If either $p > q$ or $p = 2$, then $H$ is a quantum permutation algebra. If $2 < p < q$, then $q = 1 (\emph{mod \,} p)$ and $H$ is a cocentral (non split) exact sequence $k \to k^{\mathbb Z_q \times \mathbb Z_q} \to H \to k\mathbb Z_p \to k$. In the last case, $H^*$ is a quantum permutation algebra. \end{proposition}

\begin{proof} When $p > q$, $H$ is one of the (self-dual) central abelian extensions constructed in \cite{gelaki}. On the other hand, when $p = 2$, $H$ fits into a central abelian exact sequence $k \to k^{\Gamma} \to H \to kF \to k$, where $\Gamma = \mathbb Z_q$, $F = D_q$, or $\Gamma = \mathbb Z_p$, $F = \mathbb Z_q \times \mathbb Z_q$. These extensions are classified in \cite{masuoka-18}. See \cite[Lemmas 1.3.9 and 1.3.11]{pqq}.
In the case $2 < p < q$, it follows from \cite[Subsection 1.4]{pqq} that $H$ fits into the prescribed exact sequence. The proposition follows from Theorem \ref{central}.  \end{proof}

We point out that in the case of a non split exact sequence $k \to k^{\mathbb Z_q \times \mathbb Z_q} \to H \to k\mathbb Z_p \to k$, there are examples of nontrivial Hopf algebras $H$ with no proper central Hopf subalgebra. The dual Hopf algebra $H^*$ can also be constructed as a twisting deformation of a dual group algebra. \end{example}

\begin{example} \emph{Dimension $16$.} It follows from \cite[Theorem 9.1]{kashina} that every cosemi\-simple Hopf algebra $H$ of dimension $16$ over $k$ fits into a central exact sequence $k \to
k^{\mathbb Z_2} \to H \to kF \to k$, where $F$ is group of order $8$. Therefore $H$ is a quantum permutation algebra.   \end{example}

\begin{example} \emph{Examples with irreducible characters of degree
$\leq 2$.} As another example (see \cite{qfp}), we get that the
nontrivial Hopf algebras $H = \mathcal A_n$ or $\mathcal B_n$,
studied by Masuoka in \cite{ma-contemp} are quantum permutation
algebras. Indeed, they fit into a central abelian exact sequence $k \to
k^{\mathbb Z_2} \to H \to kF \to k$, where $F$ is a dihedral
group.

It follows from \cite{BN} that if $H$ is a nontrivial semisimple
Hopf algebra and $\chi \in H^*$ is a faithful self-dual
irreducible character of degree $2$, then $H$ fits into a central
exact sequence  $k \to k^{\mathbb Z_2} \to H \to kF \to k$, where
$F$ is a polyhedral group. Therefore $H$ is also a quantum
permutation algebra in this case.

\medbreak More generally, let  $H$ be a semisimple Hopf algebra
such that its irreducible corepresentations are of
dimension $\leq 2$. Suppose in addition that $H^*$ contains no
proper central Hopf subalgebra. By \cite[Theorem 6.4]{BN} $H$ fits
into a central abelian extension $k \to k^\Gamma \to H \to kF \to
k$, with $\Gamma \neq 1$. Therefore $H$ is a quantum permutation
algebra. \end{example}

\section{Right coideal subalgebras in split extensions} Let $H = k^{\Gamma} \# kF$ be a split abelian extension.
Our aim in this section is to give some restrictions on the associated actions $\vartriangleright$ and $\vartriangleleft$, in order that $H$ contains a commutative right coideal subalgebra. In the case of right coideal subalgebras which are 'extremal' in a certain sense, we obtain conditions that correspond, roughly, to the assumptions in Theorems \ref{central} and \ref{f-abelian} (see Proposition \ref{co-pi}).
The results will be used in the next section.

\medbreak Consider the canonical projection
$\pi = (\epsilon \otimes \id): H \to kF$. Note that $\pi(e_g \# x)
= \delta_{g, 1}x$, for all $g \in \Gamma$, $x \in F$.

Then $H$ is a right $kF$-comodule algebra via $\rho = (\id \otimes
\pi) \Delta: H \to H \otimes kF$. In other words, $H$ is an
$F$-graded algebra $H = \oplus_{x \in F}H_x$, where, for all $x
\in F$,
$$H_x = \{ h \in H|\, (\id \otimes \pi)\Delta(h) = h \otimes x \} = k^\Gamma \#x.$$

Suppose $R \subseteq H$ is a right coideal subalgebra, that is,
$\Delta(R) \subseteq R \otimes H$. Then $R$ is a $kF$-subcomodule
algebra of $H$, thus it is a graded subalgebra, $R = \oplus_{x \in
F}R_x$, where $R_x = R\cap H_x$, for all $x \in F$.

\medbreak Let $\Supp R \subseteq F$ denote the \emph{support} of
$R$, that is, $\Supp R = \{ x \in F|\, R_x \neq 0\}$.

\medbreak Since $R$ is a right coideal subalgebra of $H$, then
$\pi(R) = kT$, where $T$ is a subgroup of $F$. On the other hand,
$\pi$ defines, by restriction, an epimorphism of right
$kF$-comodules $\pi: R \to kT$. In other words, $\pi: R \to kT$ is
a (surjective) map of $F$-graded spaces, with respect to the
natural grading on $kT$. Therefore
\begin{equation}\label{pirx}\pi(R_x) = \begin{cases} kx, \quad
{\rm if \,} x \in T,\\ 0, \quad {\rm
otherwise.}\end{cases}\end{equation}

Notice that, since $R_x = R \cap H_x = R \cap (k^\Gamma \# x)$,
then every nonzero element of $R_x$, $x \in \Supp R$, is of the
form $f \# x$, where $f\in k^\Gamma$ is nonzero.

Furthermore, if $x \in F$, then $x \in T$ if and only if there
exists $f \# x \in R_x$ with  $f(1) \neq 0$.

\medbreak Let $\rightharpoonup: \Gamma \times k^{\Gamma} \to
k^{\Gamma}$ denote the action by algebra automorphisms of $\Gamma$
on $k^\Gamma$ given by right translations, that is, $(s
\rightharpoonup f)(g) = f(gs)$, for all $s, g \in \Gamma$.

\begin{lemma}\label{supp-stable} \emph{(i)} For all $g \in \Gamma$, $x \in F$, we have
\begin{equation*}R_{g \triangleright x} = \{ (g \rightharpoonup f) \# (g
\triangleright x)|\; f \# x \in R_x \}. \end{equation*}

\emph{(ii)} $\Supp R$ is a $\Gamma$-stable subset of $F$
containing $T$.

\emph{(iii)} For all $x \in \Supp R$, there exists $t \in \Gamma$
such that $t \triangleright x \in T$. \end{lemma}

\begin{proof} (i). Let $x \in F$ and let $f \in k^\Gamma$, such that $f\#x \in
R_x$. Since $\Delta(R) \subseteq R \otimes H$,
\begin{align*}\Delta(f\#x) & = \sum_{g\in \Gamma} f(g) \Delta(e_g \# x)
= \sum_{s, t \in \Gamma} f(st)(e_s \# t \triangleright x) \otimes
(e_t \# x)\\ & = \sum_{t \in \Gamma} ((t \rightharpoonup f) \# t
\triangleright x) \otimes (e_t \# x) \in R \otimes H.\end{align*}

Fix $g \in \Gamma$. Evaluating the right tensorand of the last
expression in $\epsilon_F \# g \in H^*$, we get that
$(g\rightharpoonup f) \# (g \triangleright x) \in R$, and since
this is homogeneous of degree $g \triangleright x$, then
$(g\rightharpoonup f) \# (g \triangleright x) \in R_{g
\triangleright x}$. This shows that $\{ (g \rightharpoonup f) \#
(g \triangleright x)|\; f \# x \in R_x \} \subseteq R_{g
\triangleright x}$.

Since $g \in \Gamma$ was arbitrary, the other inclusion follows
from this applied to $g^{-1} \in \Gamma$. This proves (i).

\medbreak (ii). By (i), $\Supp R$ is $\Gamma$-stable. The
inclusion $T \subseteq \Supp R$ follows from \eqref{pirx}.

\medbreak (iii). Let $x \in \Supp R$ and let $f\#x \in R_x$, where
$0 \neq f \in k^\Gamma$. By (i), we have that $(t\rightharpoonup
f) \# (t \triangleright x) \in R_{t \triangleright x}$, for all $t
\in \Gamma$.

Since $f \neq 0$, there exists $s \in \Gamma$ such that $f(s) \neq
0$. Then $(s\rightharpoonup f) (1) = f(s) \neq 0$. Hence $s
\triangleright x \in T$. This proves (iii). \end{proof}

\begin{remark} It follows from Lemma \ref{supp-stable} that, for all $x \in F$, $x \in \Supp R$ if and only if $x^{-1} \in
\Supp R$. Indeed, if $s \in \Gamma$ is such that $s \triangleright
x \in T$, then $(s \triangleleft x) \triangleright x^{-1} = (s
\triangleright x)^{-1} \in T$. Thus $(s \triangleleft x)
\triangleright x^{-1} \in \Supp R$. Since $\Supp R$ is
$\Gamma$-stable, it follows that $x^{-1} \in \Supp R$. This proves
the claim.
\end{remark}

Assume in addition that $R \subseteq H$ is a \emph{commutative}
right coideal subalgebra. Then $T$ is an abelian subgroup of $F$.

\begin{proposition}\label{co-pi} Let $H = k^{\Gamma} \# kF$. Let also $R \subseteq H$ be a commutative right coideal subalgebra and let $\pi(R) = kT$,  where $T$ is an abelian subgroup of $F$.  Then:
\begin{enumerate} \item[(i)] \label{1}If $k^\Gamma \subseteq R$, then $T$ acts trivially on $\Gamma$ via $\vartriangleleft$.
\item[(ii)] \label{2} If $k^\Gamma \cap R = k1$, then $T$ is stable under the action $\vartriangleright$ of $\Gamma$.
\end{enumerate} \end{proposition}

\begin{proof} (i). Consider the $F$-gradings $H = \oplus_{x \in F}H_x$, $R = \oplus_{x \in F}R_x$, as before.
In this case, we have $k^\Gamma \subseteq R_1 \subseteq  H_1 =
k^\Gamma$. Hence $k^\Gamma = R_1 = H_1$. In particular, since
$R_x$ is an $R_1$-module under left multiplication, then $k^\Gamma
R_x \subseteq R_x$, for all $x\in F$.

Let $f\#x \in R_x$ and let $g \in \Gamma$. Then $e_g (f\# x) =
f(g) \, e_g \# x \in R_x$. Thus $R_x = k^{S} \# x$, for some
subset $S \subseteq \Gamma$. (Note that $x \in T$ if and only if
$1 \in S$.)

\medbreak Let $x \in \Supp R$ and put $R_x = k^{S} \# x$, $S
\subseteq \Gamma$, as above. Let $s \in S$. Since $e_s \in
k^\Gamma \subseteq R$, and $R$ is commutative, we have $$e_s \# x
= e_s (e_s \# x) = (e_s \# x) e_s = \delta_{s \triangleleft x, s}
e_s \# x.$$  Therefore $s \triangleleft x = s$, for all $s \in S$.


\medbreak Now suppose $x \in T$, so that $e_1\# x \in R_x$. By
Lemma \ref{supp-stable} (i), $R_{t \triangleright x} = k^{St^{-1}}
\# (t \triangleright x)$, for all $t \in \Gamma$. In particular,
$e_{t^{-1}} \# t \triangleright x = t \rightharpoonup e_1 \# t
\triangleright x \in R_{t \triangleright x}$, for all $t \in
\Gamma$.

As we have seen above, this implies that $t^{-1} \triangleleft (t
\triangleright x) = t^{-1}$, for all $t \in \Gamma$. Thus
$$1 = (t^{-1}t) \triangleleft x = \left( t^{-1} \triangleleft (t
\triangleright x)\right) (t \triangleleft x) = t^{-1} (t
\triangleleft x),$$ for all $t \in \Gamma$. This means that the
action of $x$ on $\Gamma$ is trivial. Since $x \in T$ was
arbitrary, this proves (i).

\medbreak (ii). Suppose $x \in T$ and let $f \# x \in R_x$
such that $f(1) = 1$. Since $x^{-1} \in T$, there exists $f' \#
x^{-1} \in R_{x^{-1}}$ with $f'(1) = 1$. The product $(f \# x) (f'
\# x^{-1})$ belongs to $R_{1} = k1$, by assumption.

On the other hand, we have $(f \# x) (f' \# x^{-1}) = f (x. f') \#
1$. Thus $$(f \# x) (f' \# x^{-1}) = f (x. f') \# 1 = (f (x.
f'))(1) \, 1 \# 1 = 1 \# 1,$$ since $(f (x. f'))(1) = f(1)f'(1 \triangleleft
x) = f(1)f'(1) = 1$. This shows that, for all $x \in T$, $R_x$
contains an invertible element $f\#x$ with inverse $f' \# x^{-1}
\in R_{x^{-1}}$.  In particular, $f \in k^\Gamma$ is invertible
and therefore $f(g) \neq 0$, for all $g \in \Gamma$.

\medbreak Let now $x \in \Supp R$. By Lemma \ref{supp-stable}
(iii),  $s \triangleright x \in T$, for some $s \in \Gamma$.
Hence, by the above, there is an invertible element $f \# s
\triangleright x \in R_{s \triangleright x}$ with $f(g) \neq 0$,
for all $g \in \Gamma$.

By Lemma \ref{supp-stable} (i), we have $(s^{-1}\rightharpoonup f)
\# x = (s^{-1}\rightharpoonup f) \# s^{-1}\triangleright (s
\triangleright x) \in R_x$. Moreover, $(s^{-1}\rightharpoonup f)
(1) = f(s) \neq 0$. Hence $x \in T$.

This shows that $\Supp R = T$ and then $T$ is $\Gamma$-stable, by
Lemma \ref{supp-stable} (ii). This proves (ii) and finishes
the proof of the proposition. \end{proof}

\begin{example}\label{esqpg-Julien} Let $\Gamma$ be a finite group acting by automorphisms on a finite group $F$ via $\triangleright: \Gamma \times F \to
F$, and let $H = k^\Gamma \# kF$ (so $\triangleleft$ is trivial in
this case). Thus $H$ is a central split abelian extension of
$k^{\Gamma}$.

Let $T \subset F$ be a subgroup and let
$$X(T) = {\rm Span}\{e_g \# g^{-1} \triangleright y, \ g \in \Gamma , \ y \in T\}.$$
Then $X(T)$ is a right coideal subalgebra of $H$ (containing
$k^\Gamma$) and it is commutative  if $T$ is abelian.


\medbreak  For appropriate choices of the abelian subgroup $T$,
the algebras $X(T)$ provide examples of commutative coideal
subalgebras $R \subseteq H$ with $\Supp R$ not necessarily
included in an abelian subgroup.
\end{example}

\section{Split extensions associated to the symmetric group} \label{sn} In this section we shall give examples of cosemisimple Hopf algebras $H$  which are not quantum permutation algebras.

Let $\mathbb S_n$ denote the symmetric group on $n$ symbols.
Let $H = k^{C_n} \# k\mathbb S_{n-1}$ be the split abelian extension associated to the matched pair $(C_n,
\mathbb S_{n-1})$ arising from the exact factorization $\mathbb S_{n} = \Gamma F$, where $\Gamma = C_n = \langle z \rangle \simeq \mathbb Z_n$, $z = (1 2 \dots n)$, and $F = \{ x \in \mathbb S_{n}|\, x(n) = n \} \simeq \mathbb S_{n-1}$.
(Actually,  $H$ is the unique, up to isomorphism, Hopf algebra fitting into an
exact sequence $k \to k^{\mathbb Z_n} \to H \to k\mathbb S_{n-1} \to k$ \cite[Theorem 4.1]{ma-calculations}.)

\begin{remark}\label{bgm}It follows from Theorem \ref{f-abelian} that $H^* =
k^{\mathbb S_{n-1}} \# kC_n$ is a quantum permutation algebra.
Note that, as for Drinfeld doubles, we have $D(H) \simeq D(H^*)$, and $D(H)^* \simeq
(D(\mathbb S_n)^*)^\sigma$ is a twisting deformation of $D(\mathbb S_n)^*$,  for a certain convolution invertible cocycle $\sigma :
D(\mathbb S_n)^* \otimes D(\mathbb S_n)^* \to k$ \cite{BGM}. \end{remark}

\medbreak  We quote the following fact on stabilizers, that follows from \cite[Lemma 3.2]{JM}.

\begin{lemma}\label{cn-sn}  We have $F_1 = \mathbb S_{n-1}$ and $F_{z^j} = \{ x \in \mathbb S_{n-1}|\, x(n-j) = n-j \} \simeq \mathbb S_{n-2}$, $1\leq j \leq n-1$. In particular, the only subgroup of $F = \mathbb S_{n-1}$ that acts trivially on $C_n$ is the trivial subgroup $\{ 1\}$. \end{lemma}

\begin{proof} It is shown in \cite[Lemma 3.2]{JM} that there are two orbits for the action of $\mathbb S_{n-1}$ on $C_n$, namely $\mathcal O_1 = \{ 1\}$ and $\mathcal O_z = \{ z, \dots, z^{n-1}\}$.
We have $F_1 = \mathbb S_{n-1}$ and $F_{z} = \{ x \in \mathbb S_{n-1}|\, x(n-1) = n-1 \} \simeq \mathbb S_{n-2}$.

Moreover, for each $1\leq j \leq n-1$, $z^j = z \vartriangleleft x_j$, where $x_j$ is the transposition $x_j = (n-1 \, n-j)$.
Therefore, $F_{z^j} = x_j^{-1} F_z x_j$ is the claimed subgroup of $F = \mathbb S_{n-1}$.  \end{proof}

\medbreak Suppose  that $n = p $ is a prime number.  Let $R \subseteq H$ be a commutative right coideal subalgebra, and let $\pi(R) = kT$, where $T \subseteq \mathbb S_{p-1}$ is an abelian subgroup. Since $R_1 = R \cap k^\Gamma$ is a right coideal subalgebra (hence a Hopf subalgebra) of $k^\Gamma$, then  $\dim R_1$ divides $|\Gamma| = p$. Therefore, we must have either $R \cap k^\Gamma = k^\Gamma$ or $R  \cap k^\Gamma = k1$.

By Proposition \ref{co-pi}, in the second case $T$ is
$C_p$-stable, while in the first case $T$ acts trivially on
$\Gamma = C_p$, and thus $T = 1$, in view of Lemma \ref{cn-sn}. This implies that  in this case $R \subseteq k^{C_p}$, by Lemma \ref{supp-stable} (iii).

\begin{lemma}\label{gen-abel} Let $p$ be a prime number and suppose that $H = k^{C_p} \# k\mathbb S_{p-1}$ is a quantum permutation algebra.
Then $\mathbb S_{p-1}$ is generated by its abelian $C_p$-stable subgroups. \end{lemma}

\begin{proof} The assumption implies that $H$ is generated by commutative right coideal subalgebras $R^1, \dots, R^N$, $N \geq 1$. Letting $kT^j = \pi(R^j)$, we have that the abelian subgroups $T^j$, $1\leq j \leq N$, generate $F = \mathbb S_{p-1}$.

By the above, either $T^j = 1$ or $T^j$ is $C_p$-stable. Hence $\mathbb S_{p-1}$ is actually generated by abelian $C_p$-stable subgroups, as claimed. \end{proof}

We can now state the main result of this section:

\begin{theorem}\label{s5} The cosemisimple Hopf algebras $H =
k^{C_4} \# k\mathbb S_{3}$ and $H =
k^{C_5} \# k\mathbb S_{4}$ are \emph{not} quantum permutation
algebras.

In particular, there exist finite dimensional cosemisimple Hopf algebras which are
not quantum permutation algebras. \end{theorem}

\begin{proof} Let $H = k^{C_5} \# k\mathbb S_{4}$. In this case, the action of $C_5 = \langle (12345) \rangle$ on $\mathbb S_4$ is written down
explicitly in Table 1 of \cite[pp. 15]{JM}. It turns out that  the only abelian $C_5$-stable subgroups  of $\mathbb
S_4$ are contained in the cyclic subgroup $\langle (1342)
\rangle$.  Thus they do not generate $\mathbb S_4$.
Lemma \ref{gen-abel} implies that $H$ is not a quantum permutation algebra.

\medbreak
Let now $H = k^{C_4} \# k\mathbb S_{3}$. In this case, the action of $C_4$ on $\mathbb S_3$ has three orbits:
\begin{equation}\label{orbits}\{ 1\}, \quad \{ (13) \}, \quad \{ (12), (23), (123), (132)\}. \end{equation} In particular, the only abelian subgroup of $\mathbb S_3$ which is $C_4$-stable is $\langle (13) \rangle \simeq \mathbb Z_2$.

We have in addition that $kG(H) = k^{C_4} \# k\langle (13) \rangle$, and $G(H) \simeq D_4$.

\medbreak Let $R \subseteq H$ be a commutative right coideal subalgebra. As before, consider the $\mathbb S_3$-grading $R = \oplus_{x\in \mathbb S_3} R_x$, where $R_x = R \cap (k^{C_4} \# x)$ and let $\pi(R) = kT$, where $T \subseteq \mathbb S_{3}$ is an abelian subgroup.

The subalgebra $R_1 = R \cap k^{C_4}$ is a right coideal subalgebra (hence a Hopf subalgebra) of $k^{C_4}$. Then either $R_1 = k1$ or $R_1 = k^{C_4}$, or $R_1 = k^{C_4/L}$, where $L = \{ 1, z^2\}$ is the only subgroup of order $2$ of $C_4$.

By Proposition \ref{co-pi}, the assumption $R_1 = k^{C_4}$ implies that $T$ acts trivially on
$C_4$, and thus $T = 1$, by Lemma \ref{cn-sn}. Also, if $R_1 = k1$, then $T$ is
$C_4$-stable and therefore $T \subseteq \langle (13) \rangle$.

In any of these cases, we obtain that $\Supp R \subseteq \langle (13) \rangle$, by Lemma \ref{supp-stable} (iii), and thus $R \subseteq kG(H)$.

\medbreak Suppose that there exists a commutative right coideal subalgebra $R$  such that $\Supp R \nsubseteq \langle (13) \rangle$. By Proposition \ref{co-pi} (iii), also $T \subsetneq \langle (13) \rangle$ and $\Supp R \neq \mathbb S_3$ (otherwise, the transposition $(13)$ would belong to the orbit of another cycle in $\mathbb S_3$).

By the above, $R_1 = k^{C_4/L}$ is of dimension $2$ and in view of \eqref{orbits}, \begin{equation}\label{suppR}\Supp R = \{ 1\} \cup \{ (12), (23), (123), (132)\},\end{equation}  since by Proposition \ref{co-pi} (ii), $\Supp R$ is $C_4$-stable.



\medbreak The subalgebra $R_1 = k^{C_4/L}$ of $k^{C_4}$ is spanned by the idempotents $e_L = e_1 + e_{z^2}$ and $e_{Lz} = e_z + e_{z^3}$.

For a subset $X \subseteq C_4$ and $f \in k^{C_4}$, let us denote $f_X = \sum_{x \in X}f(x)e_x$. Note that $f_L =  e_L f$, for all $f \in k^{C_4}$.

Let $x \in \Supp R$ and let  $f \in k^{C_4}$ such that $f\#x \in R_x$. Since $R$ is commutative by assumption, we have
\begin{equation}\label{fl}e_L (f \# x) = f_L \# x = (f \# x)e_L = f_{L \vartriangleleft x^{-1}} \# x, \end{equation}
that is, $f_L = f_{L \vartriangleleft x^{-1}}$. Hence, if $x \neq 1$, then $f(z^2) = 0$. Otherwise, $z^2 \vartriangleleft x = z^2$, implying that $x \in F_{z^2} = \langle (13) \rangle$ (see Lemma \ref{cn-sn}), which contradicts \eqref{suppR}.

\medbreak Now suppose that $x \in T$, $x \neq 1$. Then there exists $f \# x \in R_x$ such that $f(1) = 1$. Thus $f_L = e_1$, by the above. In particular, $e_1 \# x = e_L(f \# x) \in R_x$, and by Proposition \ref{co-pi} (i), $e_{t} \# t^{-1} \vartriangleright x \in R_{t^{-1} \vartriangleright x}$, for all $t \in C_4$.

By \eqref{fl}, this implies that  $(e_{t})_L = (e_{t})_{L \vartriangleleft (t^{-1} \vartriangleright x)^{-1}}$, for all $t \in C_4$.
In particular, taking $t = z^2$, we get that $z^2 \in L \vartriangleleft (z^2 \vartriangleright x)^{-1}$, so that $z^2 \vartriangleright x \in F_{z^2} = \langle (13) \rangle$. This contradicts again \eqref{suppR}, since $1 \neq z^2 \triangleright x \in \Supp R$.

\medbreak This shows that there can exist no commutative right coideal subalgebra $R$ with $\Supp R \nsubseteq \langle (13) \rangle$. Then we conclude that every commutative right coideal subalgebra of $H$ is contained in $kG(H)$. Therefore $H$ is not a quantum permutation algebra, by Theorem \ref{gen-cmc}.  This finishes the proof of the theorem. \end{proof}

\begin{remark} The results in Section \ref{examples-mp} ensure that any semisimple Hopf algebra of dimension less than or equal to $23$ is a quantum permutation algebra. Theorem \ref{s5} implies that this bound is optimal, since  $H = k^{C_4} \# k\mathbb S_{3}$ is a non quantum permutation algebra of dimension $24$. In particular, $24$ \emph{is the smallest possible dimension that a non quantum permutation algebra can have}. \end{remark}

\begin{remark} Observe that in the case where $H = k^{C_5} \# k\mathbb S_{4}$, the exact factorization
$\mathbb S_5 = \mathbb S_4 C_5$ restricts to an exact
factorization $\mathbb A_5 = \mathbb A_4 C_5$. We may therefore
consider the split extension $H' = k^{C_5} \# k\mathbb A_{4}$.
(Indeed, $H'$ is isomorphic to a Hopf subalgebra of $H$, by Remark
\ref{h-sub}.)

The arguments used so far apply \textit{mutatis mutandi} to this
new matched pair. Then, as before, we get that $H'$ is not a
quantum permutation algebra. The example provided by $H'$  has dimension $60$.
\end{remark}

\begin{remark} Note that if $H$ is not a quantum permutation algebra, then
the tensor product $\tilde H = H \otimes H^*$ is not a quantum
permutation algebra neither (since $H$ is a quotient of $\tilde
H$), and $\tilde H$ is \emph{self-dual}. Theorem \ref{s5}
implies that there exist self-dual cosemisimple Hopf algebras which
are not quantum permutation algebras. \end{remark}

\begin{remark} As noted before, for $H$ as in Theorem \ref{s5}, we have $D(H)^* \simeq (D(\mathbb
S_n)^*)^\sigma$ ($n = 4$ or $5$) is a twisting deformation of a quantum permutation
algebra (see Corollary \ref{dual-double} and Remark \ref{bgm}).  Since  $H$ is a quotient of $D(H)^*$, then $D(H)^*$ is not a quantum permutation algebra.

This provides an example of a twisting deformation of a quantum permutation algebra which is not a quantum permutation algebra. \end{remark}

As a consequence of Theorem \ref{s5}, and since $H^*$ is a
quantum permutation algebra, we get the following:

\begin{corollary}\label{st-dual}  The class of quantum permutation algebras in not stable neither under duality nor under Hopf algebra extensions nor under twisting deformations. \end{corollary}

\section{Quantum permutation algebras obtained from twisting}\label{twist}

We have seen in the previous section that the class of quantum permutation algebras is not stable
under twisting. We begin by giving  a stability result (Proposition \ref{suff-twist}) for twistings of quantum permutation algebras, under a technical condition on the cocycle.
Then we give (Theorem \ref{twist-abelian}) a construction of quantum permutation algebras by the twisting of certain linear algebraic groups, using the results from \cite{tga}.


The results combined together cover the known quantum permutation
algebras obtained by twisting.

\begin{proposition}\label{suff-twist}
 Let $H$ be a quantum permutation algebra generated by the coefficients of a magic matrix
 $x=(x_{ij})\in M_n(H)$. Let $\sigma : H\otimes H \longrightarrow k$ be a 2-cocycle satisfying
 $$\sigma(x_{ij},x_{il}) = \delta_{ij}\delta_{il}, \forall i,j,l$$
 Then $H^\sigma$ is a quantum permutation algebra.
\end{proposition}

\begin{proof}
Recall  that the Hopf algebra $H^\sigma$ is $H$ as a coalgebra, and
the product is defined by
$$[x] [y]= \sigma(x_{1}, y_{1})
\sigma^{-1}(x_{3}, y_{3}) [x_{2} y_{2}], \quad x,y \in H,$$ where
an element $x \in H$ is denoted $[x]$, when viewed as an element
of $H^{\sigma}$. Then
\begin{align*}
[x_{ij}][x_{il}] &= \sum_{r,s,p,q}\sigma(x_{ir}, x_{is})\sigma^{-1}(x_{pj},x_{ql})[x_{rp}x_{sq}] \\
&=\sum_{p,q} \sigma^{-1}(x_{pj},x_{ql})[x_{ip}x_{iq}] \\
&=\sum_{p}\sigma^{-1}(x_{pj},x_{pl})[x_{ip}] \\
&= \delta_{jl}[x_{ij}].
\end{align*}
since we also have $\sigma^{-1}(x_{ij},x_{il}) =
\delta_{ij}\delta_{il}, \forall i,j,l$. As $\sum_{i}[x_{ij}]=[1]=
\sum_{i}[x_{ji}]$, we conclude that $([x_{ij}])$ is a magic matrix
and hence $H^\sigma$ is a quantum permutation algebra.
\end{proof}

\begin{example}
 The twistings of $k^{\mathbb S_n}$ in \cite{bi1} are of this type.
\end{example}

Let $\Gamma$ be an abelian group and let $\sigma \in
Z^2(\Gamma,k^*)$. Recall that the character group
$\widehat{\Gamma}$ acts faithfully by automorphisms on the twisted
group algebra $k_\sigma\Gamma$ by
$$\chi. g= \chi(g)g, \ \forall \chi \in \widehat{\Gamma}, g \in \Gamma$$
So we consider $\widehat{\Gamma}$ as a subgroup of ${\rm
Aut}(k_\sigma\Gamma)$.
\begin{theorem}\label{twist-abelian}
 Let $\Gamma$ be an abelian group and let $\sigma \in Z^2(\Gamma,k^*)$. Consider a linear algebraic group
 $G$ such that $\widehat{\Gamma} \subseteq G \subseteq {\rm Aut}(k_\sigma\Gamma)$. Then $\sigma$ induces
 a 2-cocycle $\sigma'$ on $\mathcal O(G)$ such that $\mathcal O(G)^{\sigma'}$ is a quantum permutation algebra.
\end{theorem}

\begin{proof}
 The cocycle $\sigma'$ is constructed in the standard way: the inclusion
 $\widehat{\Gamma} \subset G$ induces a surjective Hopf algebra map $\mathcal O(G) \rightarrow k^{\widehat{\Gamma}}$
 which, composed with the canonical isomorphism $k^{\widehat{\Gamma}} \simeq k\Gamma$ yields a surjective Hopf algebra map
 $p: \mathcal O(G) \longrightarrow k\Gamma$. The 2-cocycle $\sigma'$ is defined by $\sigma'= \sigma(p \otimes p)$.

 Now let $A_{aut}(k_\sigma\Gamma)$ be the universal Hopf algebra coacting on $k_\sigma\Gamma$ and leaving the canonical trace invariant (see \cite{wang,tga}). Since the automorphism group ${\rm Aut}(k_\sigma\Gamma)$ preserves the canonical trace, the universal property of $A_{aut}(k_\sigma\Gamma)$ yields a Hopf algebra map
 $q: A_{aut}(k_\sigma\Gamma) \rightarrow \mathcal O(G)$. Now the composition of surjective Hopf algebra maps
 $$A_{aut}(k_\sigma\Gamma)  \overset{q} \to  \mathcal O(G) \overset{p} \to k\Gamma$$
 yields a composition of surjective Hopf algebra maps
 $$A_{aut}(k_\sigma\Gamma)^{\sigma''}  \overset{q} \longrightarrow  \mathcal O(G)^{\sigma'} \overset{p}  \longrightarrow k\Gamma,$$
 where $\sigma'' = \sigma(pq \otimes pq)$.
 We know from \cite{tga} that $A_{aut}(k_\sigma\Gamma) \simeq A_{aut}(k\Gamma)^{\sigma^{-1}(\pi \otimes \pi)}$
 where $\pi : A_{aut}(k\Gamma) \longrightarrow k\Gamma$ is the Hopf algebra map arising from the coaction of $k\Gamma$ on itself. It is then clear that $A_{aut}(k_\sigma\Gamma)^{\sigma''} \simeq A_{aut}(k\Gamma)$. The latter is a quantum permutation algebra since $k\Gamma$ is commutative, and hence we conclude that so is $\mathcal O(G)^{\sigma'}$.
\end{proof}

\begin{remark} The cocycle $\sigma'$ in Theorem \ref{twist-abelian} is the one 'lifted' from the cocycle $\sigma$ on  $\widehat{\Gamma} \simeq \Gamma$. \end{remark}

\begin{remark}
 If $k_\sigma \Gamma$ is non commutative and if the only subgroup of $\widehat{\Gamma}$ that is normal
 in $G$ is trivial, then the algebra $\mathcal O(G)^{\sigma'}$ is non commutative (see \cite{ni}).
\end{remark}

\begin{example} Let $\Gamma = \mathbb Z_2^n$, $n \geq 2$, and consider the bicharacter $\sigma$ on $\Gamma$ given by $\sigma(t_i, t_j) = -1$, if $i < j$, $\sigma(t_i, t_j) = 1$, if $i \geq j$, where $t_i$, $1\leq i \leq n$, denote the standard generators of $\Gamma$. In this case the twisted group algebra $k_{\sigma}\Gamma$ is isomorphic to the Clifford algebra $Cl_n(k) = k[t_i, 1 \leq i \leq n|\; t_i^2 = 1, \, t_it_j = -t_jt_i, \, i \neq j]$.
Since the orthogonal group $O_n(k)$ acts naturally on $Cl_n(k)$ by algebra automorphisms, we get from Theorem \ref{twist-abelian} that for any subgroup $\mathbb Z_2^n \subseteq G \subseteq O_n(k)$, the cocycle twist $\mathcal O(G)^{\sigma'}$ is a quantum permutation algebra.

When $G = O_n(k)$, we get the hyperoctahedral quantum group $O_n^{-1}(k)$ from \cite{ahn}.

\end{example}

\begin{example} Let $\Gamma = \mathbb Z_n \times \mathbb Z_n$ and let $\sigma$ be the $2$-cocycle on $\Gamma$ given by $\sigma((i, j), (t, l)) = w^{jt}$, where $w \in k^*$ is a primitive $n$-th root of unity. In this case we have $k_{\sigma}\Gamma \simeq M_n(k)$, so that $\textrm{Aut}(k_{\sigma}\Gamma) \simeq {\rm PGL}_n(k)$.

By Theorem \ref{twist-abelian}, for every linear algebraic group $G$ such that $\widehat{\Gamma} \subseteq G \subseteq {\rm PGL}_n(k)$, the cocycle twist $\mathcal O(G)^{\sigma'}$ is a quantum permutation algebra.
The twisted examples from \cite{qfp} are of this type for $n=2$.
\end{example}

\section{Quantum permutation envelope}\label{QPE}

Let $H$ be a cosemisimple Hopf algebra. Consider the
subalgebra $H_{qp} \subseteq H$ generated by the matrix coefficients
of all separable commutative (right \emph{and} left) coideal subalgebras of $H$.

\begin{lemma}\label{qpe-hopf} $H_{qp}$ is a Hopf subalgebra of $H$ containing all quantum permutation algebras $A \subseteq H$. \end{lemma}

\begin{proof} It is clear that $H_{qp}$ is a subbialgebra. Since the image  of a separable commutative right (respectively, left) coideal subalgebra under the antipode of $H$ is a separable commutative  left (respectively, right) coideal subalgebra, then we have $\mathcal S(H_{qp}) = H_{qp}$. Then $H_{qp}$ is a Hopf subalgebra.

Since every quantum permutation algebra $A \subseteq H$ is generated by matrix coefficients of some separable commutative left coideal subalgebras, then $A \subseteq H_{qp}$. This finishes the proof of the lemma. \end{proof}

\begin{definition} The Hopf subalgebra $H_{qp}$ will be called the
\emph{quantum permutation envelope} of $H$. \end{definition}

\begin{proposition}
Suppose that $H$ is finite dimensional and cosemisimple.  Then $H_{qp}$ is  the maximal quantum permutation algebra contained in $H$. Moreover, $H_{qp}$ is generated as an algebra by the matrix coefficients of all separable commutative right (or left) coideal subalgebras of $H$. \end{proposition}

\begin{proof} Note that, being finite dimensional, the subbialgebra $H'$ generated as an algebra by the matrix coefficients of all separable commutative right (or left) coideal subalgebras of $H$ is a Hopf subalgebra. This implies that $H_{qp} = H'$.

It follows from Theorem  \ref{gen-cmc} that $H_{qp}$ is a quantum permutation algebra, and it is maximal by Lemma \ref{qpe-hopf}. This proves the proposition.  \end{proof}

\begin{example} Suppose $H=\mathcal R(G)$ is the algebra of representative function on  a compact group $G$. Then the quantum permutation envelope $H_{qp} \subseteq H$
coincides with $\mathcal R(G/N)$ where $N$ is the intersection of all
closed normal subgroups of $G$ of finite index.

If $G$ is connected, the only such subgroup is $G$, so that
$H_{qp} = k1$.

On the other extreme, the condition $G/N = G$ means exactly that $G$
is \emph{residually finite}, that is, morphisms from $G$ to finite
groups separate points of $G$.

In particular, if $G$ is a  \emph{profinite group} (equivalently,
a totally discontinuous compact group \cite[I.1]{serre}), then $G$
is residually finite.
\end{example}

Regarding split abelian extensions, we have:

\begin{proposition} Let $p$ be a prime number and let $H = k^{\mathbb Z_p} \#
kF$, where $F$ is a finite group. Then $H_{qp} \subseteq H$ is the
Hopf subalgebra generated, as an algebra, by $k^{\mathbb Z_p} \# kF'$ and
$k^{\mathbb Z_p} \# kF''$, where $F'$ is the largest subgroup of $F$
acting trivially on $\mathbb Z_p$, and $F''$ is the subgroup generated by
the abelian $\mathbb Z_p$-stable subgroups of $F$. \end{proposition}

The subgroup $F'$ is the kernel of the  homomorphism $F \to
\mathbb S_{p-1}$ induced by $\triangleleft$.

\begin{proof} By Remark \ref{h-sub}, both $F'$ and $F''$ are
$\mathbb Z_p$-stable subgroups of $F$,  and the bicrossed products
$k^{\mathbb Z_p} \# kF'$ and $k^{\mathbb Z_p} \# kF''$ are Hopf subalgebras of
$H$. Since $F'$ acts trivially on $\mathbb Z_p$, then $k^{\mathbb Z_p} \# kF'$
is a central extension.

By Theorem \ref{central} and Theorem \ref{f-abelian}, $k^{\mathbb Z_p} \#
kF'$ and $k^{\mathbb Z_p} \# kF''$ are both quantum permutation algebras.
Hence the subalgebra $\tilde H \subseteq H$ generated by them is a
quantum permutation algebra, and therefore $\tilde H \subseteq
H_{qp}$. On the other hand, by Proposition \ref{co-pi}, $H_{qp}
\subseteq \tilde H$. This proves $H_{qp} = \tilde H$, as claimed.
\end{proof}

As an example, let $H = k^{C_5}\# k\mathbb S_4$ be the Hopf algebra in Theorem \ref{s5}. Then
the quantum permutation envelope of $H$ is the split extension
$H_{qp} = k^{C_5}\# k\langle (1342) \rangle$, which is a
cocommutative Hopf subalgebra (indeed, the action of $C_5$ on $\langle (1342) \rangle$ is trivial, as follows from \cite{JM}), with $\dim H_{qp} = 20$.

It follows from the definition of the actions that $H_{qp} \simeq kG$, where $G \simeq \mathbb F_5 \rtimes \mathbb F_5^*$.

\medbreak We have seen in the proof of Theorem \ref{s5} that if $H = k^{C_4}\# k\mathbb S_3$, then every commutative right coideal subalgebra of $H$ is contained in $kG(H)$. Hence, in this case, $H_{qp} = kG(H) \simeq kD_4$.

\bibliographystyle{amsalpha}

\end{document}